\documentclass{llncs}
\usepackage{graphicx}
\usepackage{wrapfig}
\usepackage[hang]{subfigure}
\usepackage{amsmath,amssymb,latexsym,graphicx,color,verbatim}
\usepackage{url} 

\title{Geometric deformations of sodalite frameworks}
\author{Ciprian Borcea\inst{1} \and Ileana Streinu\inst{2}}
\institute{
Department of Mathematics, Rider University, Lawrenceville, NJ 08648, USA
\and
Computer Science Department, Smith College, Northampton, MA 01063, USA}

\begin{document}
\maketitle

\begin{abstract}
In mathematical crystallography and computational materials science, it is important to infer flexibility properties of framework materials from their geometric representation. We study combinatorial, geometric and kinematic properties for frameworks modeled on sodalite.
\end{abstract} 
 
\medskip  \noindent
{\bf Keywords:} periodic framework, sodalite, geometric flexibility.

\section{Introduction}

In this paper we study the deformation space of a periodic framework modeled after the ideal sodalite.   

The sodalite framework is the prototype of numerous crystal structures. In fact, more than 900 crystals of this family have been identified \cite{FB}. 
This structure is remarkable not only for crystallography and materials science. It has direct connections with tiling and sphere packing problems and has natural
generalizations in arbitrary dimension \cite{O'K}.

A one-parameter geometric deformation of sodalite was first observed by Pauling \cite{P} in 1930, and used as a phase transition model \cite{M,T,Dep}.  

We show that, besides this classical `Pauling tilt scenario', the deformation space of sodalite includes a six-dimensional component, which can be described in fairly intuitive terms. 

Our study relies on the mathematical foundations of a deformation theory for periodic frameworks introduced in \cite{BS1}. In general, the deformation space is a semi-algebraic set made of real solutions of a finite system of polynomial equations. For applications, such as displacive phase transitions in crystalline materials, it is highly desirable to obtain explicit and  detailed descriptions of the deformation space \cite{D,BS2,BS3}.

 We start with the presentation of the three-dimensional structure of the ideal sodalite framework in a realization  with congruent regular tetrahedra and maximal crystallographic symmetry \cite{M,KS,BS3}. We express the {\em periodic graph} underlying the framework via 6-rings of tetrahedra with marked periods.  This formulation leads naturally into considerations of symmetry preservation in deformations, with central symmetry at the forefront. The existence of a six-dimensional deformation component follows naturally from this formulation. We also investigate dihedral symmetry and use it to explain the one-dimensional `classical Pauling tilt'.

\section{A placement with maximal symmetry}

\begin{figure}[h]
 \centering
  \subfigure[]{\includegraphics[width=0.4\textwidth]{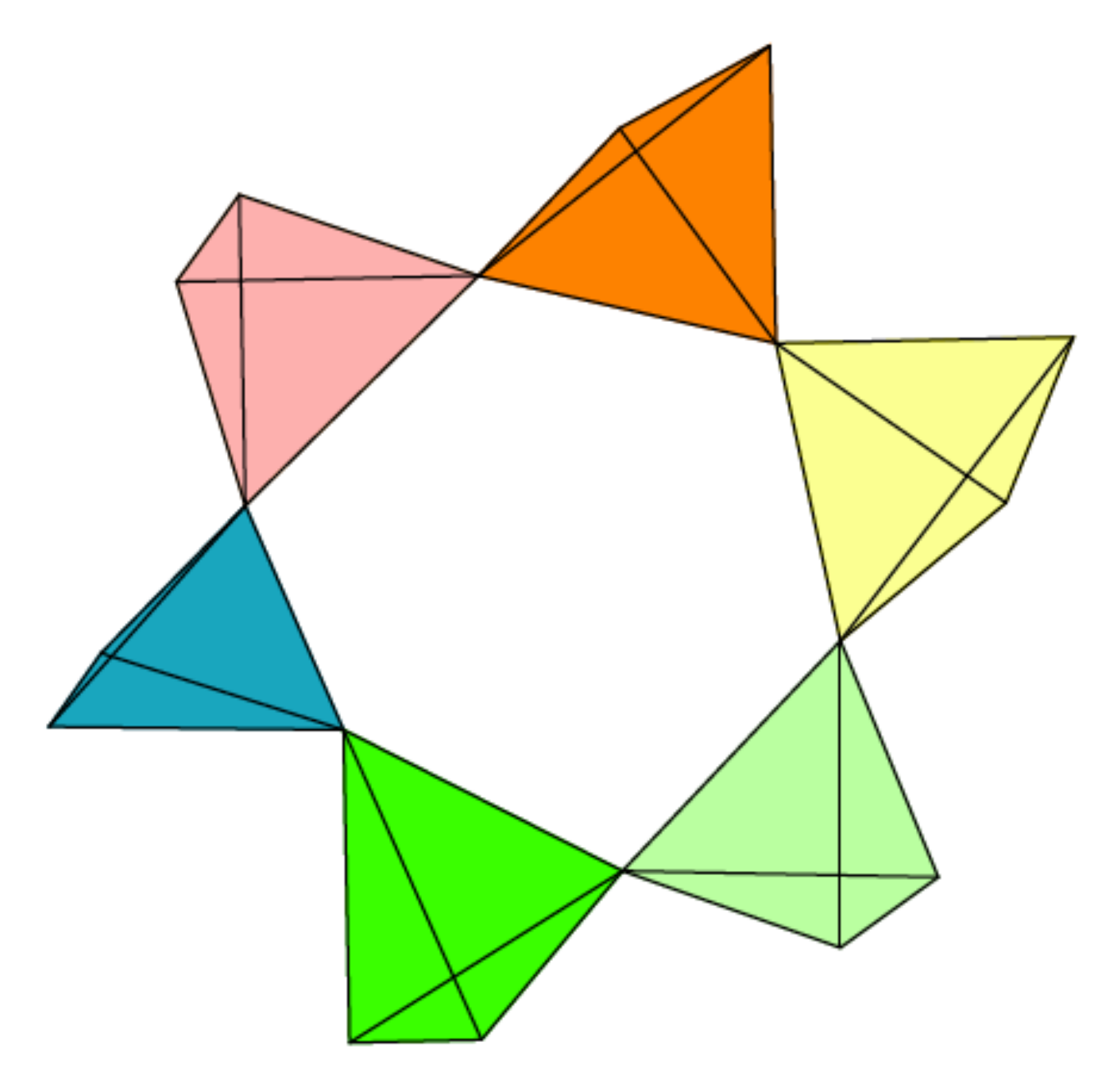}}
  \subfigure[]{\includegraphics[width=0.4\textwidth]{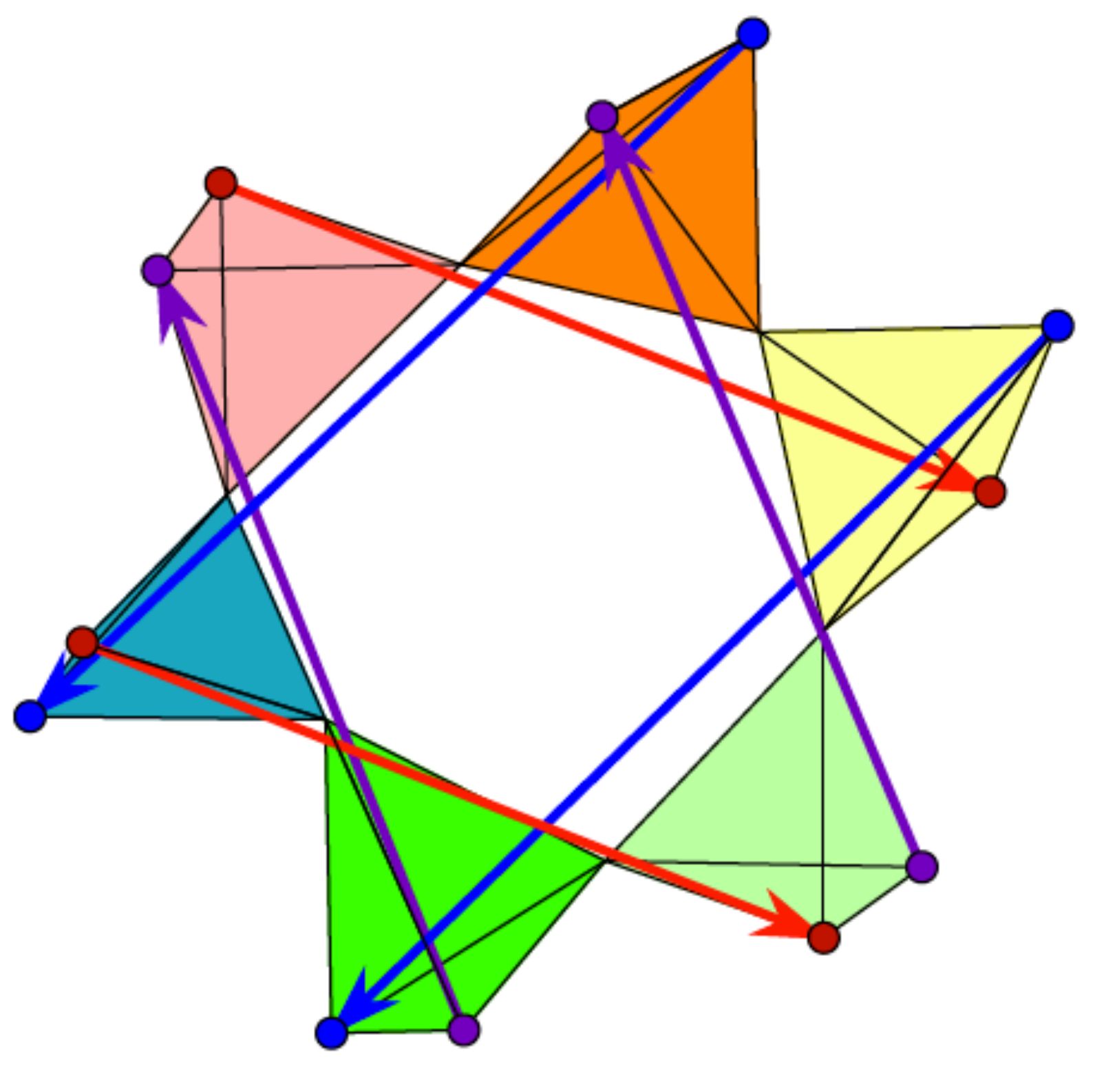}}\\
  \subfigure[]{\includegraphics[width=0.45\textwidth]{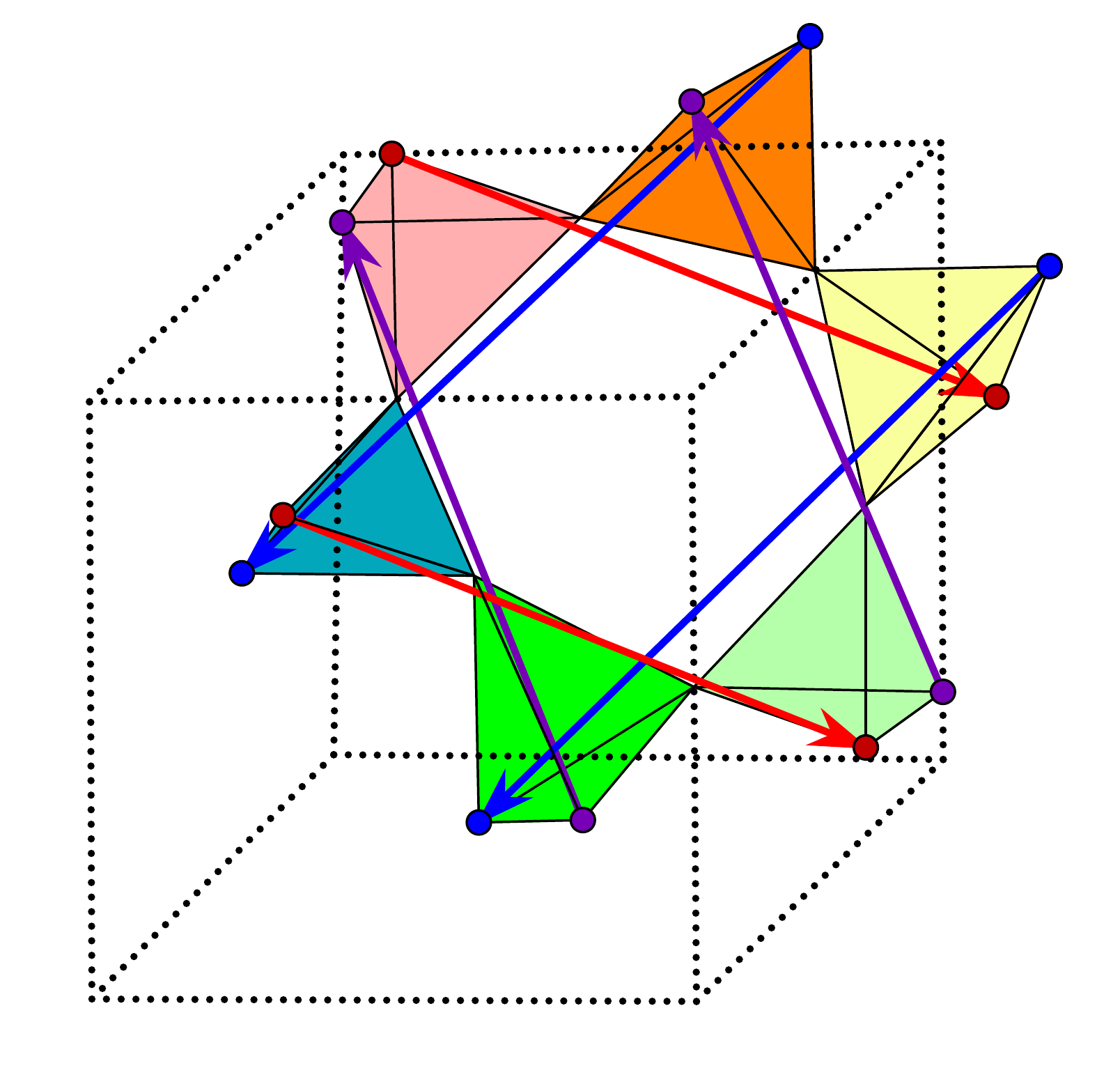}}
  \subfigure[]{\includegraphics[width=0.45\textwidth]{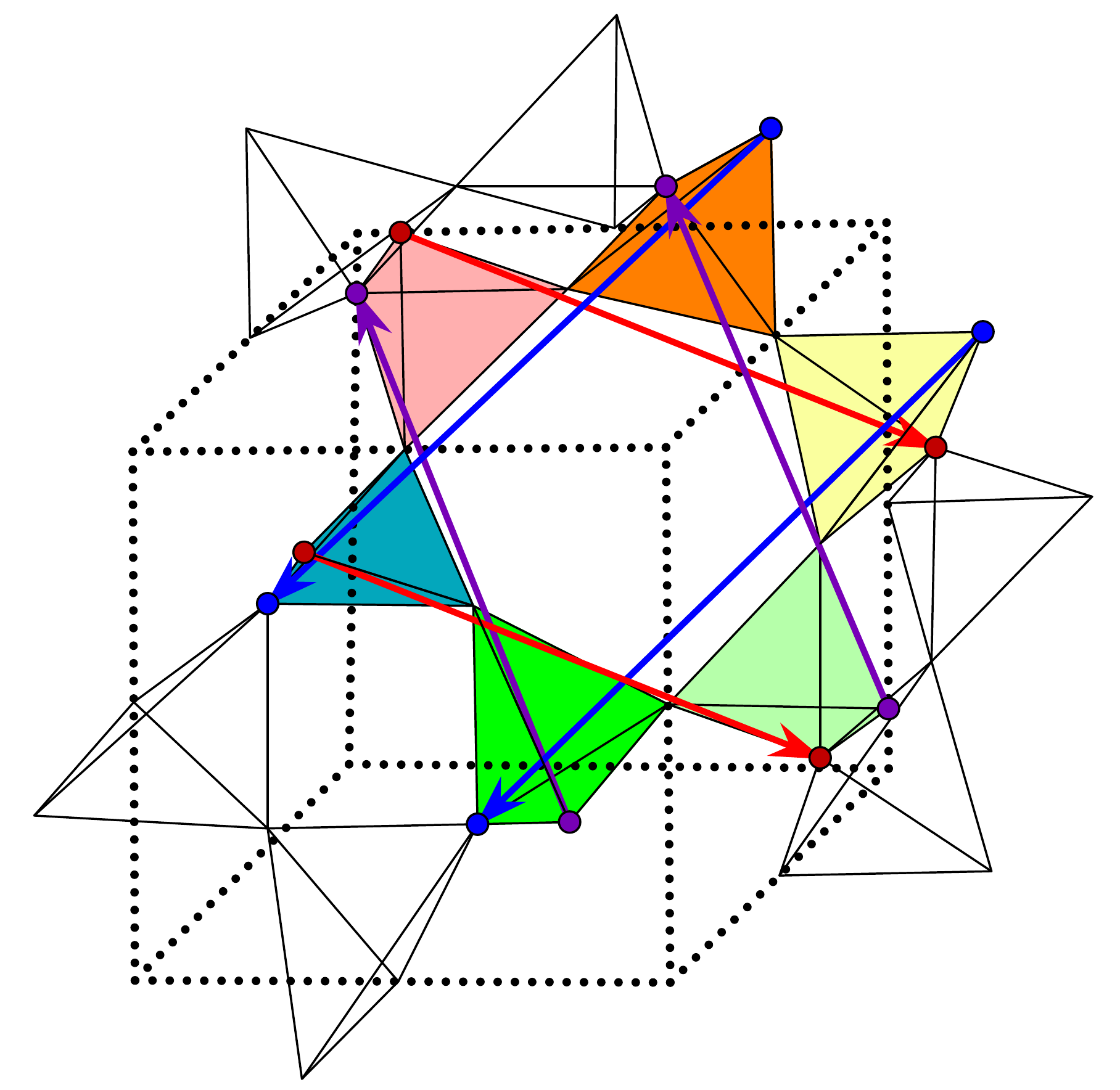}}
 \caption{An ideal sodalite placement with cubical symmetry. A 6-ring of this framework is highlighted, with three pairs of generators of the periodicity lattice. In each pair the vectors must be equal. All tetrahedra are congruent and regular.}
 \label{FigCubicS}
\end{figure}

We use Figure~\ref{FigCubicS} to illustrate the essential aspects of the initial placement of the sodalite framework: a 6-ring of regular tetrahedra (a) is placed in a specific manner relative to a cube of side length 2 centered at the origin (c). Three pairs of parallel vectors between specific vertices of the ring (b) yield the generators of the periodicity lattice.

\newpage
\noindent
Specifically, we consider the origin of the Cartesian system of coordinates at the center of the cube depicted in Fig.\ref{FigCubicS} (c), with standard basis vectors $e_1, e_2, e_3$  reaching to the centers of the frontal, right side and top face respectively. Thus, the edge length for the cube is 2. The edge length for the regular tetrahedra of the framework is determined by the condition that each tetrahedron has a pair of opposite (and hence orthogonal) edges parallel to a pair of standard directions $e_i,e_j$. This requires an edge lengths of $2a$, with $a=\sqrt{2}-1$.

\medskip \noindent
{\bf The $4\times 6=24$ tetrahedra around the cube.} We give the coordinates for the four
vertices $v_0, ..., v_3$ of the tetrahedron touching the cube at $v_0=e_1+e_3$. They are:

\begin{equation}\label{tetra}
v_0=(1,0,1), \ v_1=(1,0,2\sqrt{2}-1),\ 
\end{equation}
$$ v_2=(\sqrt{2}-1,\sqrt{2}-1,\sqrt{2}), \ v_3=(\sqrt{2}-1, 1-\sqrt{2},\sqrt{2}) $$ 

\medskip \noindent
The coordinates for the vertices of the other tetrahedra are obtained by symmetry. The symmetries of the cube are represented by permutations and sign changes of the three coordinates and form a group of order $3!\cdot 2^3=48$. This group has a double-transitive
action on the 24 framework tetrahedra around the cube. 

\medskip \noindent
{\bf Periodicity lattice.}\ The translational symmetries of the whole framework are generated by the three (pairs of equal) vectors shown in Figure~\ref{FigCubicS}(b). In coordinates:

\begin{equation}\label{periods}
\lambda_1=\sqrt{2}(1,-1,-1),\ \lambda_2=\sqrt{2}(-1,1,-1),\ \lambda_3=\sqrt{2}(-1,-1,1)
\end{equation}

\noindent
with the natural labeling suggested by the order three rotational symmetry around the 
diagonal direction $e_1+e_2+e_3$. We note that the periods

\begin{equation}\label{index2}
 \lambda_2+\lambda_3=-2\sqrt{2}e_1,\  \lambda_3+\lambda_1=-2\sqrt{2}e_2,\   \lambda_1+\lambda_2=-2\sqrt{2}e_3
\end{equation}

\noindent
are {\em mutually orthogonal} and generate a {\em sublattice of index 2} in the full lattice of periods generated by (\ref{periods}).

\medskip
\noindent
{\bf The sodalite cage and the Kelvin polyhedron.} The 24 framework tetrahedra wrapped around a cube form a so-called {\em sodalite cage}. We note that the barycenter of the tetrahedron in 
(\ref{tetra}) has coordinates $b=(\frac{1}{\sqrt{2}},0,\sqrt{2})$ and the barycenters of the 6-ring of tetrahedra highlighted in Figure~\ref{FigCubicS} form the vertices of a planar regular hexagon of edge one, centered at $c=\frac{1}{\sqrt{2}}(1,1,1)$.

\begin{figure}[h]
 \centering
 {\includegraphics[width=0.5\textwidth]{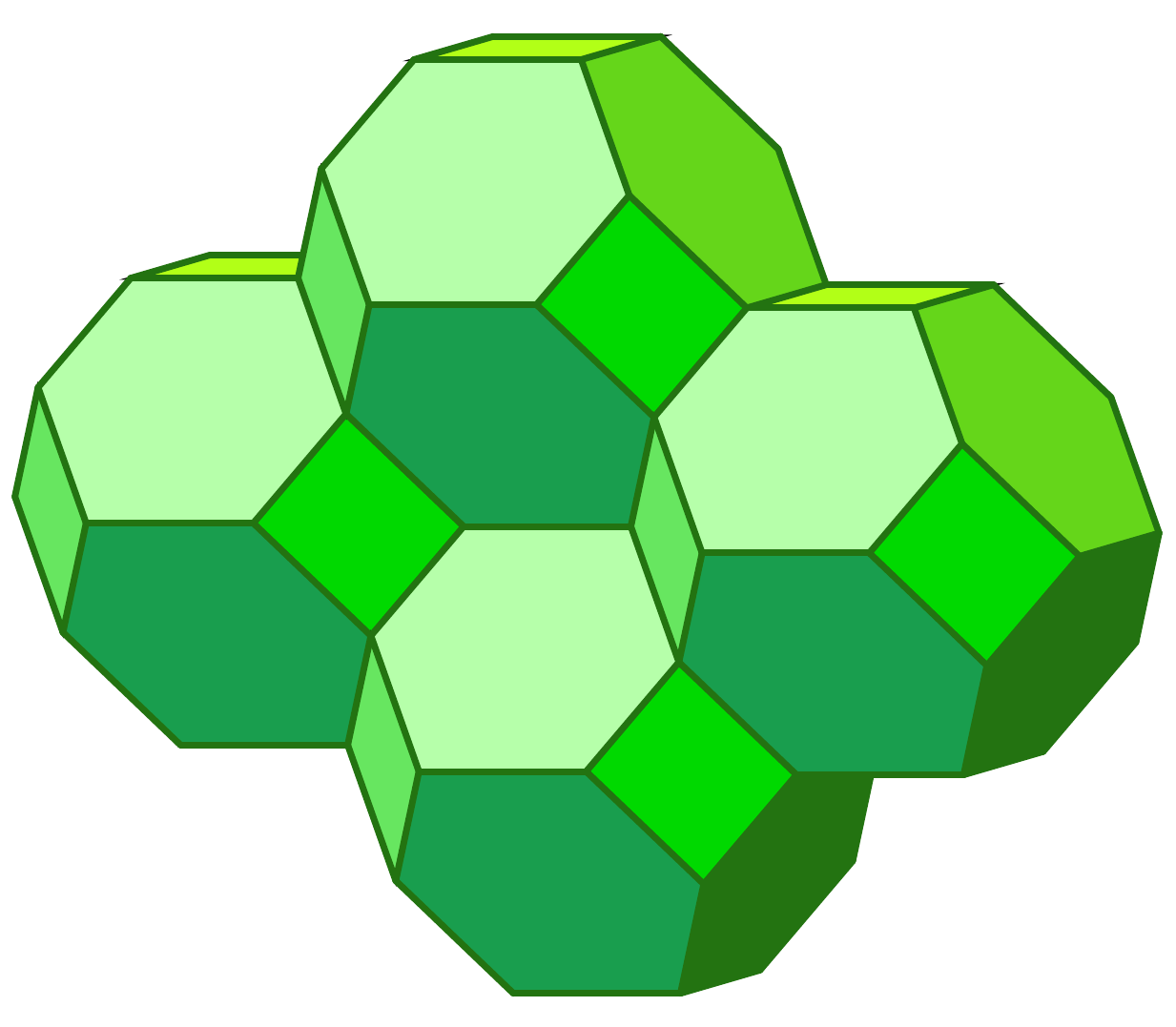}}
 \caption{Kelvin polyhedra tiling space.}
 \label{FigKelvinTiling}
\end{figure}

Thus, the 24 barycenters of the sodalite cage tetrahedra give the vertices of a centrally symmetric polyhedron with 8 hexagonal faces and 6 square faces. Under translation by the 
periodicity lattice $\Lambda$
generated by (\ref{periods}), this polyhedron {\em tiles} the three dimensional space $R^3$, as illustrated in Figure~\ref{FigKelvinTiling}. Translations identifying opposite faces are periods and generate the periodicity lattice $\Lambda$. 
The polyhedron is the {\em Voronoi cell} at the origin for the lattice $\Lambda$.
It may be seen as a truncated regular octahedron (with edge length 3 at our chosen scale).
It is also called the Kelvin polyhedron, due to a famous conjecture formulated by Kelvin \cite{Th} and disproved by a counter-example of Weaire and Phelan \cite{WP}. The name {\em permutohedra} is also used for a family which includes this polyhedron.

\vspace{-14pt}
\begin{figure}[h]
 \centering
 {\includegraphics[width=0.6\textwidth]{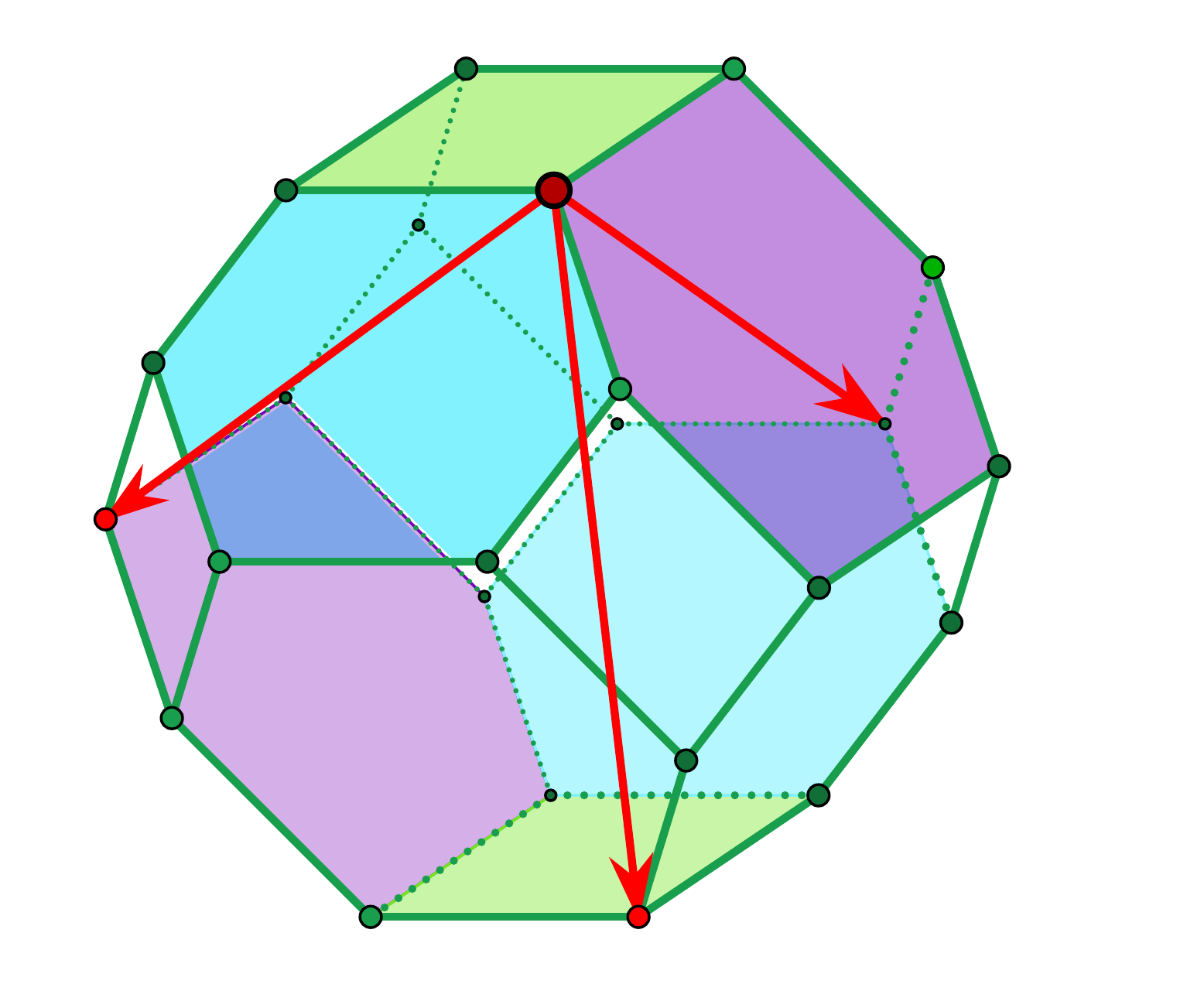}}
 {\includegraphics[width=0.3\textwidth]{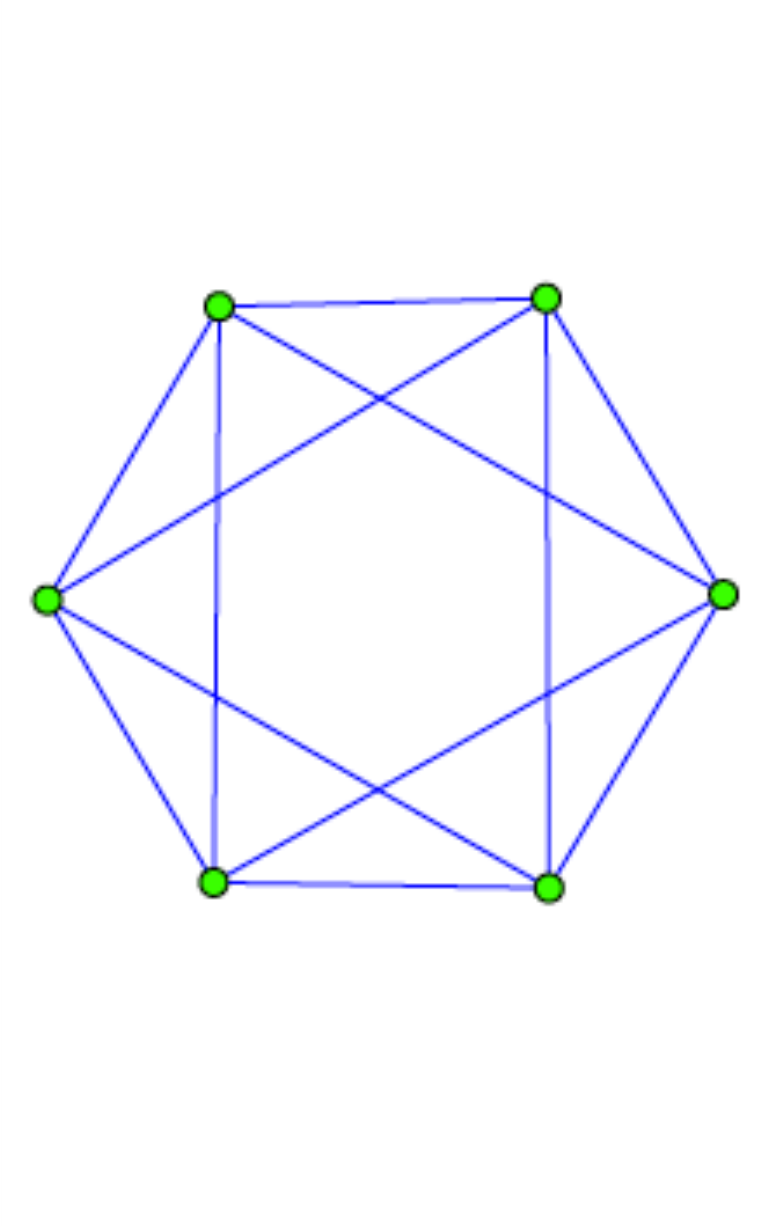}}
 \vspace{-14pt}
 \caption{(Left) Periodicity identifications in the Kelvin polytope. Opposite faces are identified by translation. (Right) The quotient graph of the Kelvin polytope skeleton.}
 \label{FigSFactorization}
\end{figure}
\vspace{-14pt}

\medskip \noindent
Figure~\ref{FigSFactorization} illustrates the result of periodicity identifications in the vertex-and-edge skeleton of a sodalite cage. The quotient graph pattern is also encoded in 
a 6-ring (with marked periods), as discussed below.

\medskip
\noindent
{\bf The 6-rings.} \ A 6-ring of tetrahedra, as in Figure~\ref{FigCubicS}, contains all the information needed for generating the whole periodic framework. It contains vertex representatives for all the vertex orbits and edge representatives for all the edge orbits under periodicity. The 6 pairs of vertices which have to be further identified by periodicity provide the generators of the periodicity lattice $\Lambda$. It is important to note that (up to sign) we have exactly three generators, since the 6 periods, as free
vectors, come as three pairs of equal vectors. This is a consequence of the 
{\em central symmetry} of the 6-ring. Indeed, $c=\frac{1}{\sqrt{2}}(1,1,1)$ is actually the center of symmetry for the 6-ring from Figure~\ref{FigCubicS}.

\vspace{-14pt}
\begin{figure}[h]
 \centering
 {\includegraphics[width=0.75\textwidth]{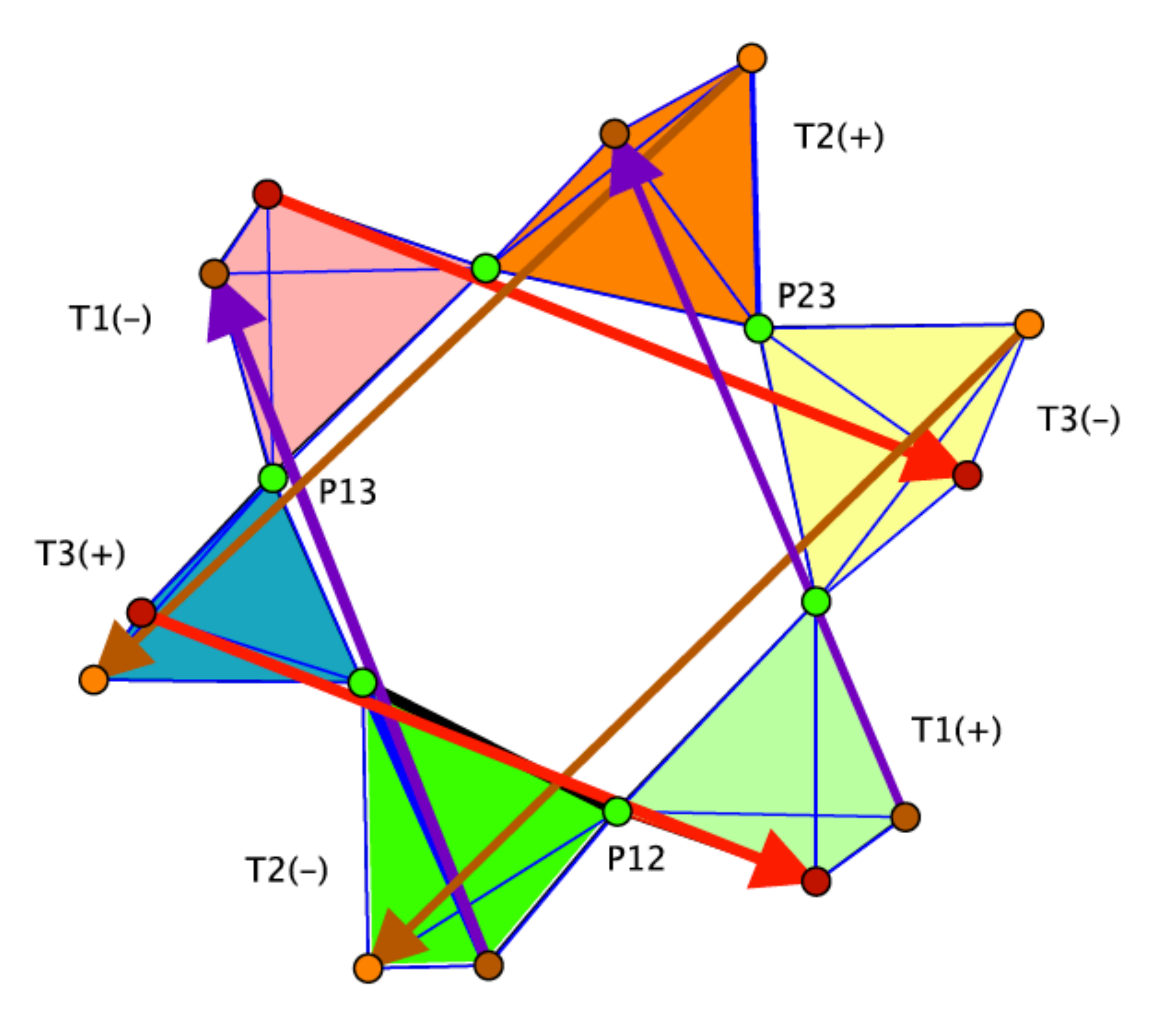}}
 \caption{ A 6-ring of the sodalite framework with three pairs of generators of the periodicity lattice.  In each pair the vectors must be equal.}
 \label{FigS}
\end{figure}

\medskip
\noindent
In order to have a convenient representation of the symmetries of the cube (given by permutations and sign changes of the coordinates) we {\bf label} the tetrahedra
in the 6-ring as follows: $T_1^-, T_3^+,T_2^-,T_1^+,T_3^-,T_2^+$, with $T_3^+,T_2^-$
on the frontal face of the depicted cube. Then coordinate sign changes will correspond to
(upper) sign changes and permutations of coordinates will correspond to permutations of
lower indices and multiplicative signature effect on upper signs. For example, the transposition of the first two coordinates corresponds to the order two product of transpositions
$(T_1^-, T_2^+)(T_3^-,T_3^+)(T_2^-,T_1^+)$. Figure~\ref{FigS} summarizes this description.

\medskip \noindent
We denote $P_{13}=T_1^{-}\cap T_3^{+}$, $Q_{23}=T_3^{+}\cap T_2^{-}$, with similar labeling around the {\bf spatial hexagon} $P_{13}Q_{23}P_{12}Q_{13}P_{23}Q_{12}$.  Note that, in the initial placement, the triangles $P_{12}P_{23}P_{13}$ and $Q_{12}Q_{23}Q_{13}$ are centrally symmetric with respect to $c=\frac{1}{\sqrt{2}}(1,1,1)$, but are not in the same plane. If we ignore the conditions on the three pairs of period vectors, the 6-ring, as a finite linkage has twelve degrees of freedom: the spatial hexagon has a six-dimensional deformation space and for a fixed configuration of the hexagon, each tetrahedron may rotate around the corresponding edge.

\section{Deformations}

We use the previously described placement for a single 6-ring to investigate the periodic deformations allowed by the ideal sodalite structure. 
Specifically, we show below that the deformation space has a six-dimensional component related to the preservation of central symmetry for the 6-ring. This ample geometrical flexibility supports the experimentally observed versatility of sodalite \cite{Dep}.

\vspace{-20pt}
\begin{figure}[h]
 \centering
 {\includegraphics[width=0.65\textwidth]{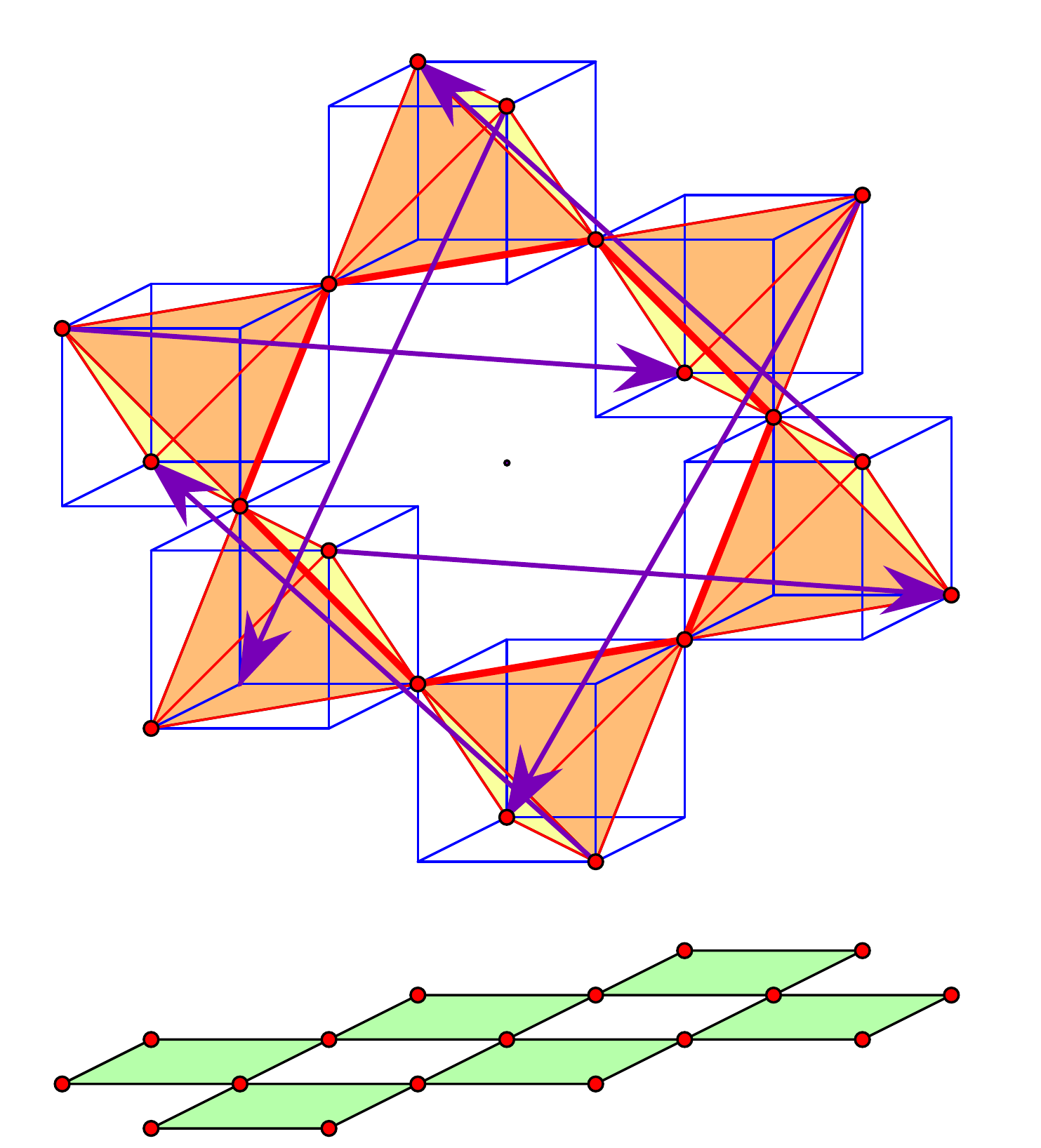}}
 \caption{ A 6-ring in a periodic deformation of the ideal sodalite framework, with three pairs of generators of the periodicity lattice. The ring has central symmetry but not periodicity preserving reflections. Cubes and their projections serve only suggestive purposes for spatial positioning. }
 \label{FigSCS}
\end{figure}

\noindent
{\bf Deformations preserving central symmetry.} This  six-dimensional component allows a very direct and intuitive description. Since we want to preserve the central symmetry of the 6-ring, indicating the configuration of three linked tetrahedra, say $T_1^{-}, T_2^{+}, T_3^{-}$,  will be enough. By equivalence under rigid motions we may assume $T_2^{+}$ fixed. Then, the positions of $T_1^{-}$ and $T_3^{-}$ can be parametrized by $SO(3)\times SO(3)$. The center of symmetry will be at the midpoint of $P_{13}Q_{13}$. The completion of the 6-ring by central symmetry will have {\em  ipso facto} equal vectors in the three pairs of periods. Thus, this component of the deformation space 
is parametrized by $SO(3)\times SO(3)$ minus the subvariety where the three generators become dependent. Topologically, the rotation group $SO(3)$ is the projective space $P_3(R)$, hence our six-dimensional component is parametrized by an open and dense subset of $P_3(R)\times P_3(R)$. One configuration is illustrated in Figure~\ref{FigSCS} and shows that reflection symmetries of the initial 6-ring are lost.

\medskip
\noindent
{\bf Deformations preserving dihedral symmetry $D_3$.} We now investigate deformations which preserve the symmetries induced by  permutations of coordinates in the initial placement. This group of order six may be conceived as a dihedral group $D_3$ generated by reflections in three planes with a common axis and forming dihedral angles of $\pi/3$ or $2\pi/3$. When restricting our attention to our chosen 6-ring, the group is simply transitive on the six tetrahedra and their barycenters form the vertices of a planar regular hexagon. 

\medskip \noindent
Suppose we want to preserve the reflection symmetry given by the transposition of the first two coordinates in the initial placement. The reflection plane goes through $P_{12}$ and $Q_{12}$ and the two marked periods
from $T_2^{-}$ to $T_1^{-}$ and respectively $T_1{+}$ to $T_2^{+}$ are parallel to this plane and one is the reflection of the other.

\medskip \noindent
Thus, in order to retain this feature when deforming the 6-ring, we look only at the three linked tetrahedra $T_1^{-}, T_3^{+}, T_2^{-}$ and adopt as reflection plane the plane through $P_{12}Q_{12}$ which runs parallel to the period from $T_2^{-}$ to $T_1^{-}$. When we complete the 6-ring by reflection in this plane, we have to satisfy only one of the remaining periodicity
constraints, since the other one will then be fulfilled as well by reflection.

\vspace{-20pt}
\begin{figure}[h]
 \centering
 {\includegraphics[width=0.65\textwidth]{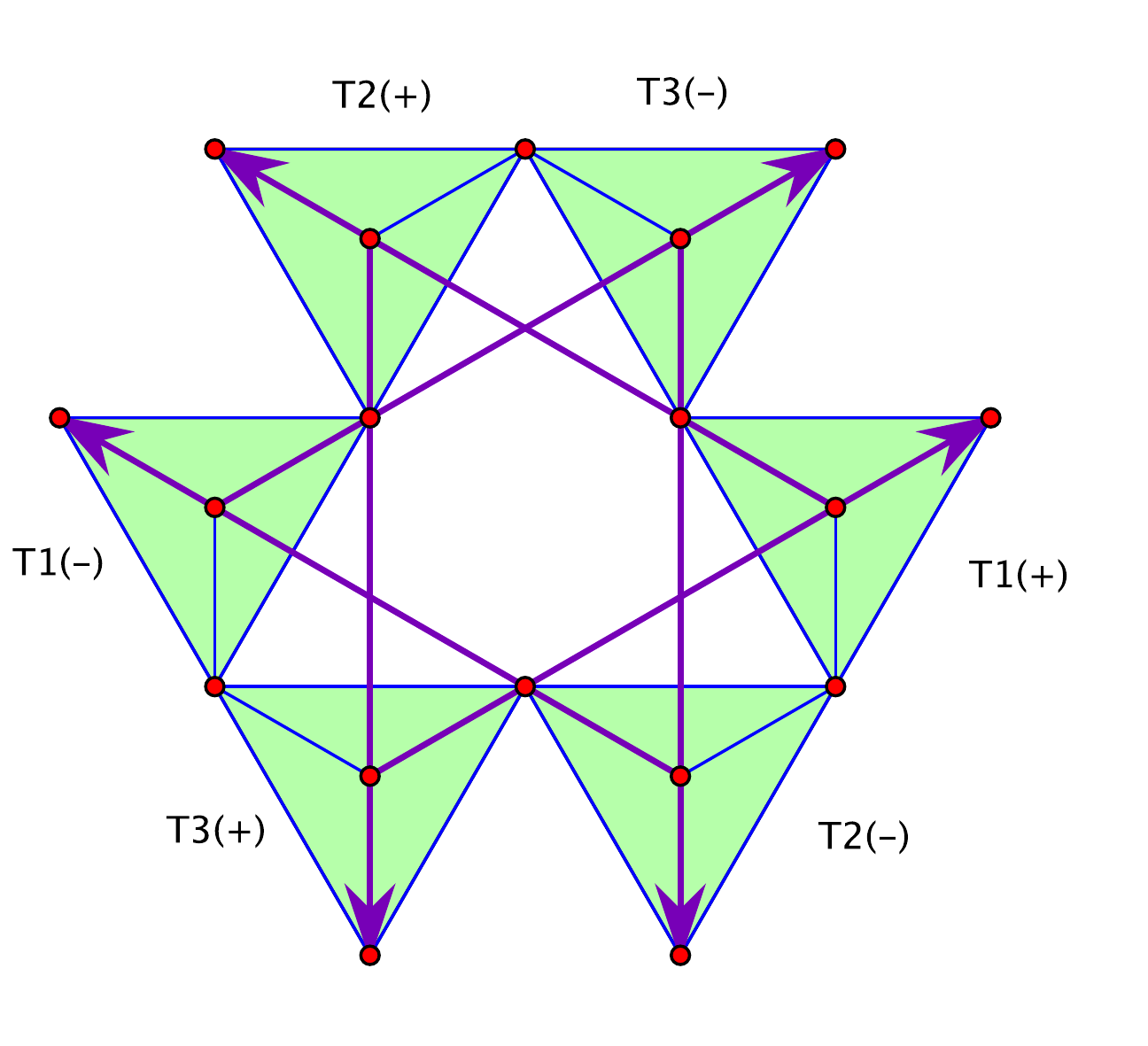}}
 \caption{ A  tetrahedrite 6-ring resulting from a periodic deformation of the ideal sodalite framework. The regular tetrahedra are pointing towards the viewer.  The 6-ring has $D_3$ dihedral symmetry but no central symmetry.}
 \label{FigT}
\end{figure}

\medskip \noindent
Similar considerations apply for the other two transpositions of coordinates. In Figure~\ref{FigT} we show a deformation which has preserved
all three reflections, but central symmetry for the 6-ring has been lost. The illustrated structure is that of another mineral called {\em tetrahedrite}. According to \cite{BF}, this relationship between sodalite and tetrahedrite was noticed only at a later stage in the mineralogy literature, by A.S. Povarennykh. We now compare deformations which `break' the central symmetry of the 6-ring and those which maintain it.

\begin{figure}[h]
 \centering
 {\includegraphics[width=0.65\textwidth]{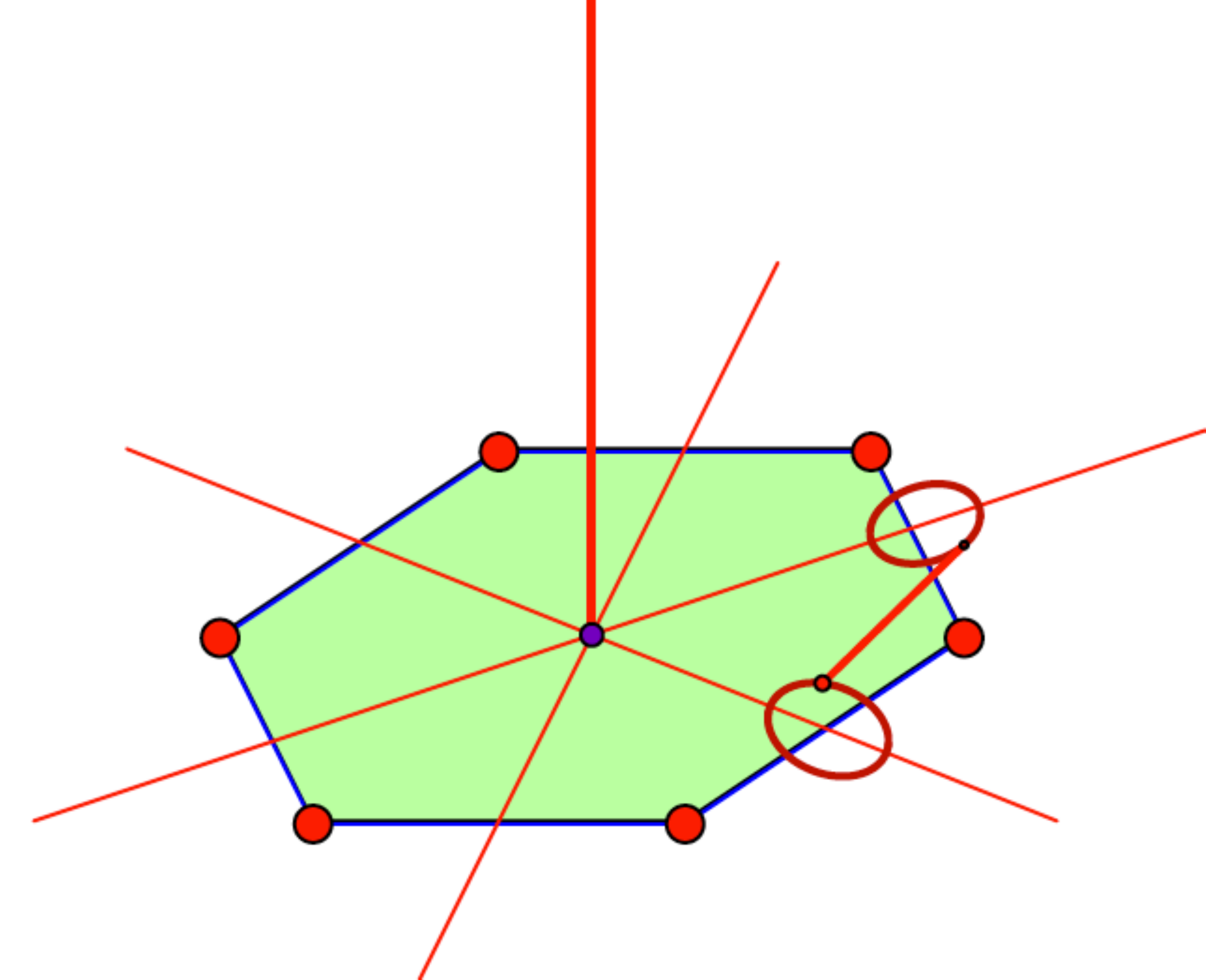}}
 \caption{ The regular hexagon of barycenters in a 6-ring with $D_3$ dihedral symmetry. Circumscribed spheres for linked tetrahedra meet in circles on the reflection planes. A possible edge position is shown.}
 \label{FigD}
\end{figure}

\noindent
Given the observed fact that barycenters must form a planar regular hexagon, we may start with the configuration shown in Figure~\ref{FigD}. The connecting vertices of the tetrahedra in the 6-ring must be in the reflection planes, indicated by their common axis and intersections with the plane of barycenters. More precisely, they belong to circles of intersection of circumscribed spheres of consecutive tetrahedra. Since one edge determines by symmetry the entire 6-ring, we see that we have a one-parameter family of possible configurations for the 6-ring with given hexagon of barycenters. When we ask for fulfillment of the periodicity conditions, the selection of the edge is subject to an additional condition expressing the parallelism of one of the periods to the corresponding reflection plane. 

\medskip
These considerations show that deformations preserving a dihedral $D_3$ symmetry consist of
curves. One type breaks the central symmetry of the 6-ring: this is the `Pauling tilt scenario' \cite{P,T,M,HG,Dep}. As mentioned above, the tetrahedrite structure can be reached via this classical tilt. What is distinctive in this scenario is that sodalite cages maintain their shape,
although decreasing in size relative to the initial case.

\begin{figure}[h]
 \centering
 {\includegraphics[width=0.50\textwidth]{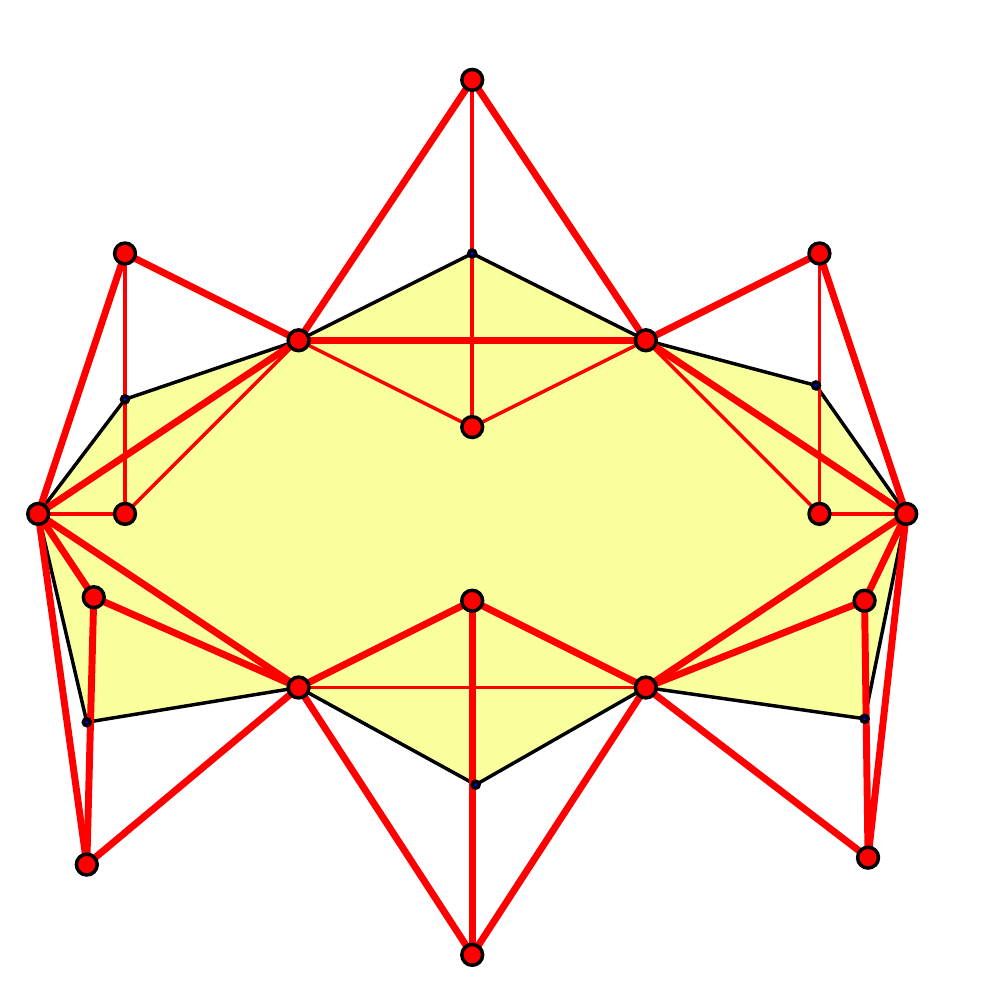}}
 \caption{ A deformation result when preserving both $D_3$ and central symmetry.
The highlighted plane is the perpendicular bisecting plane of all distant edges.}
 \label{FigF}
\end{figure}

\medskip \noindent
However, central symmetry may be maintained as follows. Since we have now double transitivity on the six tetrahedra,  all diameters in the hexagon of barycenters must go
through midpoints of opposite edges in the corresponding tetrahedra. This condition is satisfied when we choose the edge discussed above to be met perpendicularly and in the middle
by the appropriate diameter. A deformation result of this procedure is illustrated in Figure~\ref{FigF}.

\vspace{-10pt}

\section{Conclusion}

The sodalite framework belongs to the class of tectosimplicial periodic structures made of vertex sharing simplices. In dimension three, a basic count of infinitesimal deformations yields at least three infinitesimal degrees of freedom for this type of frameworks \cite{BS1}, Theorem 4.2. However, if more than three infinitesimal degrees of freedom are present, as is the case for the ideal sodalite framework, infinitesimal considerations are not sufficient for obtaining actual local deformations. In this paper, we have relied on direct and intuitive geometric features to show that for ideal sodalite, in addition to the one-parameter deformation known as the `Pauling tilt', 
there is a six-dimensional deformation component. It is conceivable that certain crystalline materials with distorted sodalite cages may be more closely related to sodalite on this account.

\end{document}